\documentclass[12pt]{imsart}
\usepackage[OT1]{fontenc}
\usepackage[utf8]{inputenc}
\usepackage[colorlinks]{hyperref}
\usepackage{amsmath,amssymb,amsthm}
\usepackage{mathrsfs}
\usepackage[english]{babel}
\usepackage{graphicx}
\usepackage{color}
\usepackage{indentfirst}
\usepackage{multirow}
\usepackage{float}
\usepackage{subfig}
\usepackage{enumerate} 
\allowdisplaybreaks

\usepackage[authoryear]{natbib}
\usepackage{a4wide}

\startlocaldefs
\numberwithin{equation}{section}
\theoremstyle{definition}
\newtheorem{theorem}{Theorem}[section]

\newtheorem{lemma}{Lemma}[section]

\newtheorem{proposition}{Proposition}[section]
\newtheorem{remark}{Remark}[section]
\endlocaldefs

\begin{document}

\newcommand\E{\mathbb{E}}
\newcommand\Var{\mathop{\text{Var}}}
\newcommand\cov{\mathbb{Cov}}
\newcommand\tr{\mathop{\text{tr}}}
\newcommand{\Sum}[2]{\sum_{\scriptstyle #1 \atop \scriptstyle #2}}
\newcommand{\rom}[1]{\uppercase\expandafter{\romannumeral #1}}
\newcommand{\g}[2]{g_{(#1)}(#2)}
\newcommand{\f}[2]{f_{(#1)}(#2)}
\newcommand{\summ}[2]{\mathop{\sum}_{#1}^{#2}}
\newcommand\re[1]{{\color{red}#1}}
\newcommand\bl[1]{{\color{blue}#1}}

\newcommand\veps{\varepsilon}

%

\begin{frontmatter}
  \title{On  singular values distribution of a
    large auto-covariance matrix in the ultra-dimensional regime}
  \runtitle{Singular values distribution of a ultra-large auto-covariance matrix}

  \thankstext{T2}{The research of J. Yao is partly supported by GRF
    Grant HKU 705413P.}
  \begin{aug}
    \author{\fnms{Qinwen} \snm{Wang}\ead[label=e1]{wqw8813@gmail.com}}
   \and
    \author{\fnms{Jianfeng} \snm{Yao}\thanksref{T2}\ead[label=e3]{jeffyao@hku.hk}}

    \runauthor{Q. Wang and J. Yao}

    \affiliation{Zhejiang  University and The University of Hong Kong}

    \address{Qinwen Wang  \\
      Department of Mathematics\\
      Zhejiang University \\
      \printead{e1}
    }

   \address{Jianfeng Yao \\
      Department of Statistics and Actuarial Science\\
      The University of Hong Kong\\
      Pokfulam, \quad
      Hong Kong \\
      \printead{e3}
    }
  \end{aug}

  \begin{abstract}
    Let $(\veps_{t})_{t>0}$ be a sequence of independent   real random
    vectors of $p$-dimension   and let $X_T=
    \sum_{t=s+1}^{s+T}\veps_t\veps^T_{t-s}/T$ be the lag-$s$ ($s$ is a
    fixed positive integer) auto-covariance matrix of $\veps_t$.
     This paper   investigates the limiting behavior of
    the singular values  of $X_T$ under the so-called {\em ultra-dimensional regime} where
    $p\to\infty$ and  $T\to\infty$ in a related way such that
    $p/T\to 0$. First, we show that the singular value  distribution
    of  $X_T$ after a suitable normalization  converges to a nonrandom
    limit $G$ (quarter law) under the forth-moment condition. Second, we
    establish the convergence of its largest singular value to the right edge
    of $G$. Both  results  are derived using the moment method.
  \end{abstract}


\end{frontmatter}

\section{Introduction}
Let $s$ be a   fixed positive integer and
$(\veps_{t})_{1\le t \le T+s}$  a sequence of independent  real
random vectors, where $\veps_t=(\veps_{it})_{1\leq i \leq p}$
has independent coordinates satisfying
$\E \veps_{it}=0$ and $\E \veps_{it}^2=1$.
Consider  the so-called  {\em lag-$s$  sample autocovariance matrix}
of $(\veps_t)$ defined as
\begin{align}
  \label{eq:acv}
  X_T=\frac1T
  \sum_{t=s+1}^{s+T}\veps_t\veps^T_{t-s}~.
\end{align}
Motivated by their application in  high-dimensional statistical
analysis where the dimensions $p$ and $T$ are assumed large (tending
to infinity), spectral analysis of such
sample autocovariance matrices have attracted much attention in recent
literature in random matrix theory.
For example, perturbation theory on the matrix $X_T$ has been carried
out  in  \citet{LY012} and \citet{Li14} for estimating the number of factors
in a large dimensional factor model of type
\begin{align}
y_t=\Lambda f_t+\veps_t+\mu~,
\end{align}
where $\{y_t\}$ is a  $p$-dimensional sequence observed at time $t$,
$\{f_t\}$ a sequence of $m$-dimensional ``latent factor'' ($m\ll p$)
uncorrelated with the error process $\{\veps_t\}$ and $\mu \in
\mathbb{R}^p$ is the general mean.
Since $X_T$ is not symmetric, its spectral distribution is
given by the set of  its singular  values which are by definition the
square roots of
positive  eigenvalues of
\begin{align}\label{at}
  A_T:=X_TX_T^T~.
\end{align}
To our best knowledge, all the existing results on $X_T$ (or $A_T$)
are found under what we will refer as
the {\em Mar\v{c}enko-Pastur regime},  or simply the
{\em MP regime}, where
\begin{align}
  p\to\infty, \quad T\to\infty \quad \text{ and} \quad
  p/T\to c>0~.
\end{align}
For example,  \citet{jin} derives the limit of the eigenvalue
distributions (ESD) of
the symmetrized
auto-covariance matrix $\frac12 ( X_T + X_T^T)$; and
\citet{wang} establishes  the  exact separation property
of the ESD which also implies  the convergence of its extreme
eigenvalues.
For the singular value distribution of $X_T$,
the limit (LSD) has been established in \citet{li} using the method of
Stieltjes transform and
in  \citet{WY14} using the moment method. The latter paper also  establishes
 the almost sure convergence of the largest singular value
of $X_T$ to the right edge of the LSD, thanks to the
moment method. Related results are also proposed in
\citet{liu} where the sequence $(\veps_t)$ is replaced by a more
general time series.

In this paper, we investigate the same questions as in  \citet{WY14}
but under a different asymptotic regime,
the so-called {\em  ultra-dimensional regime}
where
\begin{align}
  \label{eq:regime}
  p\to\infty, \quad  T\to\infty \quad \text{and}
  \quad p/T\to 0.
\end{align}
It is naturally  expected that the limit under this regime will be
much  different than under  the  MP regime above.
The findings of the paper confirm this difference by providing a new
limit of the singular value distribution of $X_T$ under the
ultra-dimensional regime.

In a related paper \citet{WangPaul14}, the authors also adopted the
ultra-dimensional regime
to derive the LSD  for a large class of separable sample covariance
matrices.
However, the autocovariance matrix $X_T$ considered in this paper
is very different of these separable sample covariance
matrices.

Recalling the definition of $A_T$ in \eqref{at}, we have
\begin{align*}
 A_T(i,j)= \frac{1}{T^2}\sum_{l=1}^p \sum_{m=1}^T\sum_{n=1}^T
 \veps_{i\,m+s}\veps_{lm}\veps_{j\,n+s}\veps_{ln}~.
\end{align*}
It follows  by simple calculations that
\begin{align*}
  \E A_T(i,j)=
  \begin{cases}
    0, & \quad  i\neq j~,\\
    p/T, & \quad i=j~,
  \end{cases}
\end{align*}
and  for $i\neq j$,
\begin{align*}
  \Var  A_T(i,j)= \E  A^2_T(i,j)=\frac{p}{T^2}~.
\end{align*}
The row sum of the variances $ \Var  A_T(i,j)$ is thus of order
$p^2/T^2$.
Therefore, in order to have the spectrum of $A_T$ be of constant order when $p/T \rightarrow 0$, we should normalise it as
\begin{align}\label{targeta}
  A:= \frac{A_T}{\sqrt{p^2/T^2}} =\frac T pX_TX_T^T~.
\end{align}

The main results of the paper are as follows. First in
Section~\ref{lsd}, we derive the almost sure  limit of the
singular value distribution of $\displaystyle \sqrt{\frac{T}p}X_T$ under the ultra-dimensional
regime and assuming that the fourth
moment of the entries $\{\veps_{it}\}$ are uniformly bounded.
This limit (LSD) simply equals to the image measure of the semi-circle law
on $[-2,2]$ by the absolute value transformation $x\mapsto |x|$.
Next in Section~\ref{largest}, we establish the almost sure
convergence of the
largest singular value  of  $\displaystyle \sqrt{\frac{T}p}X_T$
to 2 assuming that the entries
$\{\veps_{it}\}$ has a uniformly bounded moment of order $4+\nu$ for
some $\nu>0$. Both results are derived using the moment method.
Some technical details on the traditional
truncation and renormalisation steps are postponed to the appendixes.

\section{Limiting spectral distribution by the  moment method}\label{lsd}
In this section, we show that when $p/T \rightarrow 0$, the ESD of the singular values of $\sqrt{\frac T p}X_T$ tends to a nonrandom limit, which is linked to the well known semi-circle law.

\begin{theorem}\label{singularlsd}
  Suppose the following conditions hold:
  \begin{itemize}
  \item[(a).]
    $(\veps_{t})_t$ is a sequence of independent $p$-dimensional real
    valued random vectors with independent entries $\veps_{it}$,
    $1\le i\le p$,  satisfying
    \begin{align}\label{four}
      \E (\veps_{it})=0,\quad\E\veps_{it}^2=1,\quad \sup_{it}\E (\veps_{it}^4)<\infty~.
    \end{align}
  \item[(b).]  Both $p$ and $T$ tend to infinity in a related way such that $p/T \rightarrow 0$.
  \end{itemize}
  Then, with probability one, the empirical distribution of the singular values of $\sqrt{\frac Tp}X_T$ tends to the quarter law $G$ with  density function
  \begin{align}\label{densityfunctionsingular}
    g(x)=\frac{1}{\pi} \sqrt{4-x^2}~, \quad\quad 0<x\leq 2~.
  \end{align}
\end{theorem}

\begin{remark}
  Recall that the quarter law $G$ is the image measure of the
  semi-circle law by the absolute value transformation.
  It is also worth noticing that if there were no lag, i.e. $s=0$,
  the matrix $X_T$ would be a standard sample covariance matrix; and
  in this case the spectral distribution of  $\displaystyle \sqrt{\frac{T}p}(X_T-I_p)$
  would converge to the semi-circle law, see \citet{BaiYin88}. The
  case of a auto-covariance matrix $X_T$ with a positive lag $s>0$ is
  then very different.
\end{remark}

 Since the singular values of $\sqrt{\frac Tp}X_T$ are the square
 roots of the eigenvalues of $\frac T p X_TX^T_T$, in the remaining of
 this paper, we focus  on  the limiting behaviours of the eigenvalues
 of $\frac T p X_TX^T_T$.  These properties can  then be transferred  to
 the singular values of $\sqrt {\frac T p}X_T$ by the square-root
 transformation $x\mapsto\sqrt x$.

 \begin{theorem}\label{th1}
   Under the same conditions as in Theorem \ref{singularlsd},
   with probability one, the empirical spectral distribution $F^{A}$ of the matrix $A$ in \eqref{targeta} tends to a limiting distribution $F$, which is the image measure of the semi-circle law on $[-2,2]$ by the square transformation. In particular, its  $k$-th moment is:
   \begin{align}\label{mka}
     m_k=\frac 1k \left(\begin{array}{c}
         2k\\
         k-1
       \end{array}\right)~,
   \end{align}
   and its Stieltjes transform $s(z)$ and density function $f(x)$ are given by
   \begin{align}\label{stieltjes}
     s(z)=-\frac 12+\sqrt{\frac 14-\frac 1z}~, \quad\quad z \notin (0, 4]~,
   \end{align}
   and
   \begin{align}\label{densityfunction}
     f(x)=\frac{1}{\pi} \sqrt{\frac 1 x-\frac 14}~, \quad\quad 0<x\leq 4~,
   \end{align}
   respectively.
 \end{theorem}

\begin{remark}
  The $k$-th moment  in \eqref{mka} is exactly the  $2k$-th moment of
  the LSD of a standard Wigner matrix, which is also the number of
  Dyck paths of length $2k$ (for the definition of Dyck paths, we
  refer to \citet{Tao12}).
  Notice also that the density function $f$ is unbounded at the origin.
\end{remark}

The remaining of the section is devoted to the proof of
Theorem \ref{th1} using the moment method.
The $k$-th moment of the ESD $F^A$  of $A$ is
\begin{align}\label{expresss}
  m_k(A)&=\frac 1p \tr A^k=\sum_{{\bf i}=1}^T\sum_{{\bf j}=1}^p\frac{1}{p^{k+1}T^{k}}\veps_{j_1\,i_1}\veps_{j_1\,i_2}\veps_{j_2\,s+i_2}\veps_{j_2\,s+i_3}\cdots\nonumber\\
  &\quad\quad\quad\quad\quad\quad \veps_{j_{2k-1}\,i_{2k-1}}\veps_{j_{2k-1}\,i_{2k}}\veps_{j_{2k}\,s+i_{2k}}\veps_{j_{2k}\,s+i_1}~.
\end{align}
Here, the indexes in ${\bf i}=(i_1, \cdots, i_{2k})$ run over $1,2,\cdots,T$ and the indexes in ${\bf j}=(j_1, \cdots, j_{2k})$ run over $1,2,\cdots,p$.

The core of the proof is to establish
the following two assertions:
\begin{align*}
  &\text{(I)}.\quad\E m_k(A)\rightarrow m_k=\displaystyle\frac 1k \left(\begin{array}{c}
      2k\\
      k-1
    \end{array}\right)~,~~k\ge 0;\\
  &\text{(II)}. \quad\sum_{p=1}^{\infty}\Var (m_k(A))<\infty~.
\end{align*}
This is given in the Subsections \ref{ssec1}, \ref{ssec2}
and \ref{ssec3} below. It follows from these assertions
that almost surely,
$m_k(A) \to m_k$ for all $k\ge 0$.
Since the limiting moment sequence $(m_k)$ clearly satisfies the
Carleman's condition, i.e.
$\sum_{k>0} m_{2k}^{-1/(2k)}=\infty$,
we deduce that almost surely, the sequence of ESDs $F^A$
weakly converges to a probability measure $F$ whose moments are
exactly $(m_k)$. Next, notice that
$  m_k$ is exactly the number of Dyck paths of length $2k$
\citep{Tao12}, which is
also the $2k$-th moment of the semi-circle law with support $[-2,2]$,
it follows that the LSD $F$ equals to the image of the semi-circle law
by the square  transformation $x \rightarrow x^2$. The formula in
\eqref{stieltjes} and \eqref{densityfunction} are thus easily derived
and the proof of Theorem~\ref{th1} is complete.

\subsection{Preliminary steps and some graph concepts}\label{ssec1}

We now introduce the proofs for Assertions (I) and (II).
First we show that with a uniformly bounded fourth order moment,
the variables $\{\veps_{it}\}$ can be truncated at rate $\eta T^{1/4}$
for some vanishing sequence $\eta=\eta(T)$.
This is justified in Appendix~\ref{j1}.
After these truncation, centralisation and rescaling steps,
we may assume in all the following that
\begin{align}\label{jus1}
  \E(\veps_{ij})=0, ~\E\veps_{ij}^2=1,~|\veps_{ij}|\leq \eta T^{1/4}~,
\end{align}
where $\eta$ is chosen such that $\eta\rightarrow 0$ but $\eta T^{1/4}\rightarrow \infty$.

Now we introduce some basic concepts for graphs associated to the big
sum in \eqref{expresss}.
Let
\begin{eqnarray}\label{atss}
  \begin{array}{c}
    \psi (e_1, \cdots, e_m):=\text{number of distinct entities among}~ e_1, \cdots, e_m~,\nonumber\\
    {\bf i}:=(i_1, \cdots, i_{2k}),\quad {\bf j}:=(j_1, \cdots, j_{2k}),\nonumber\\
    1\leq i_a \leq T,\quad  1\leq j_b \leq p,\quad a,b=1, \cdots, 2k,\nonumber\\
    A(t,s):=\{({ \bf i},{\bf j}): \psi({\bf i})=t, \psi({\bf j})=s\}~.
  \end{array}
\end{eqnarray}
Define $Q(i,j)$ as the multigraph as follows:
Let $I$-line, $J$-line be two parallel lines, plot $i_1, \cdots, i_{2k}$ on the $I$-line, $j_1, \cdots, j_{2k}$ on the $J$-line, called the {\em $I$-vertexes} and {\em $J$-vertexes}, respectively. Draw $k$ down edges from $i_{2u-1}$ to $j_{2u-1}$, $k$ down edges from $i_{2u}+s$ to $j_{2u}$, $k$ up edges from $j_{2u-1}$ to $i_{2u}$, $k$ up edges from $j_{2u}$ to $i_{2u+1}+s$ (all these up and down edges are called {\em vertical edges}) and $k$ {\em horizontal edges} from $i_{2u}$ to $i_{2u}+s$, $k$ horizontal edges from $i_{2u-1}+s$ to $i_{2u-1}$ (with the convention that $i_{2k+1}=i_1$), where all the $u$'s are in the region: $1\leq u\leq k$. An example of the multi-graph $Q(i,j)$ with $k=3$ is presented in the following Figure \ref{cha}.
\begin{figure}[h!]
  \centering
  \includegraphics[width=13cm]{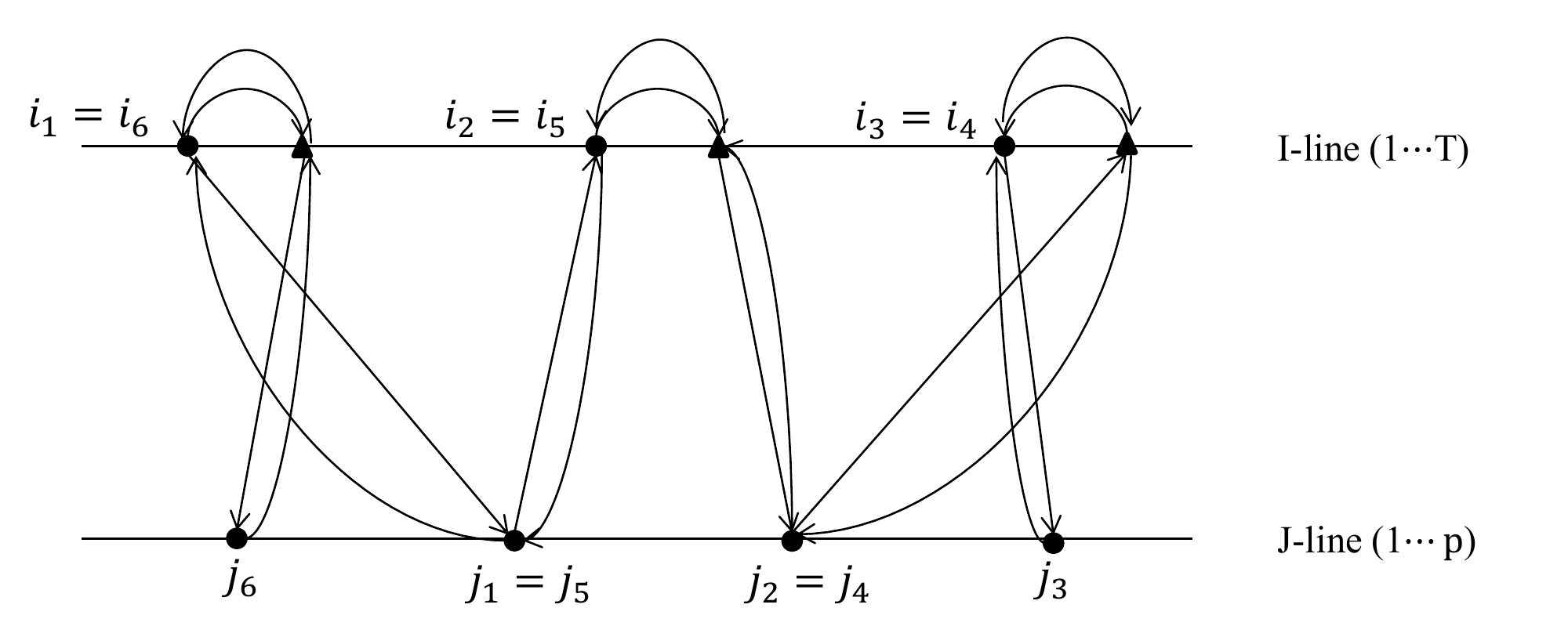}
  \caption{\small An example of the multigraph $Q(i,j)$ with $k=3$.}\label{cha}
\end{figure}

In the graph $Q(i,j)$, once a $I$-vertex $i_l$ is fixed, so is $i_l+s$. For this reason, we glue all the $I$-vertexes which are connected through horizon edges and denote the resulting graph as $M(A(t,s))$, where  $A(t,s)$ is the index set that has $t$ distinct $I$-vertexes and $s$ distinct $J$-vertexes. An example of $M(A(3,4))$ that corresponds  to the $Q(i,j)$ in Figure \ref{cha} is presented in the following Figure \ref{ma}.
\begin{figure}[h!]
  \centering
  \includegraphics[width=13cm]{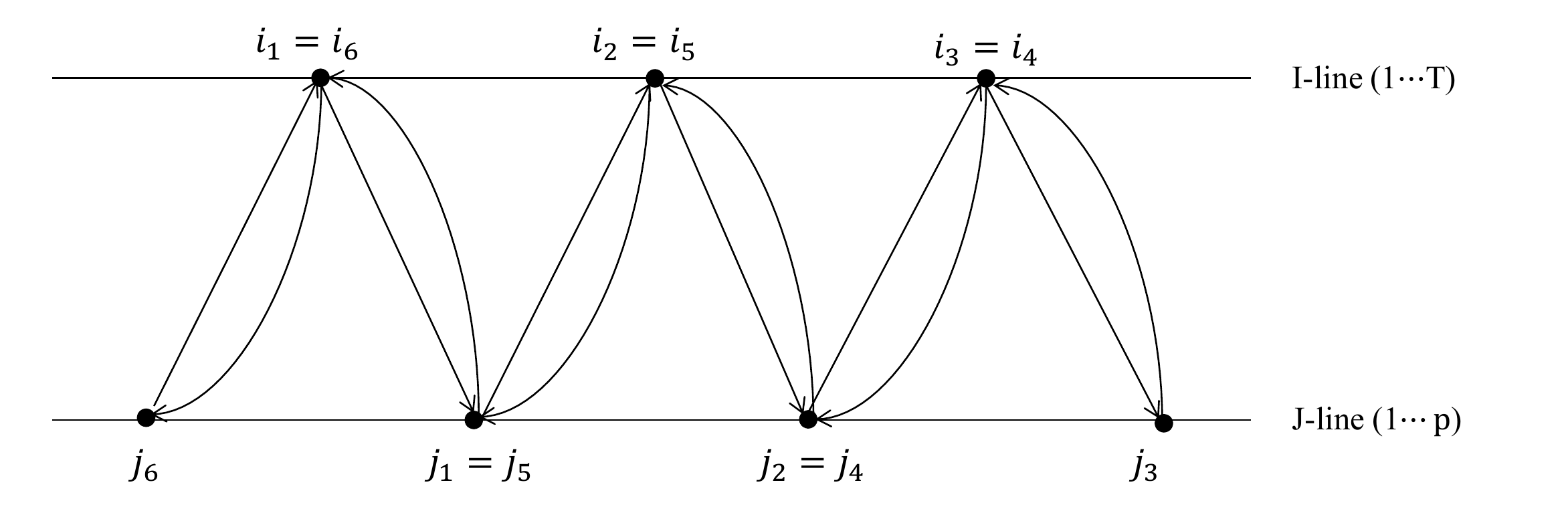}
  \caption{\small An example of $M(A(3,4))$ that corresponds  to the $Q(i,j)$ in Figure \ref{cha}.}\label{ma}
\end{figure}

\subsection{Proof of Assertion (I)}\label{ssec2}

Recall
the expression of $m_k(A)$ in \eqref{expresss}, we have
\begin{align}\label{emka}
  \E m_k(A)&= \sum_{{\bf i}=1}^T\sum_{{\bf j}=1}^p\frac{1}{p^{k+1}T^{k}}\E\big[\veps_{j_1\,i_1}\veps_{j_1\,i_2}\veps_{j_2\,s+i_2}\veps_{j_2\,s+i_3}\veps_{j_3\,i_3}\veps_{j_3\,i_4}\veps_{j_4\,s+i_4}\veps_{j_4\,s+i_5}\nonumber\\
  &\quad\quad\quad\quad\quad\quad\cdots \veps_{j_{2k-1}\,i_{2k-1}}\veps_{j_{2k-1}\,i_{2k}}\veps_{j_{2k}\,s+i_{2k}}\veps_{j_{2k}\,s+i_1}\big]\nonumber\\
  &=\sum_{t,s}\frac{1}{p^{k+1}T^k}\sum_{M(A(t,s))}p(p-1)\cdots (p-s+1)T(T-1)\cdots (T-t+1)\nonumber\\
  &\quad\quad\quad\quad\quad\quad\cdot\E\big[\veps_{j_1\,i_1}\veps_{j_1\,i_2}\veps_{j_2\,s+i_2}\veps_{j_2\,s+i_3}
  \cdots \veps_{j_{2k-1}\,i_{2k}}\veps_{j_{2k}\,s+i_{2k}}\veps_{j_{2k}\,s+i_1}\big]\nonumber\\
  &:=\sum_{t,s}S(t,s)~,
\end{align}
where
\begin{align}\label{sts}
  S(t,s)&=\frac{1}{p^{k+1}T^k}\sum_{M(A(t,s))}p(p-1)\cdots (p-s+1)T(T-1)\cdots (T-t+1)\nonumber\\
  &\quad\quad\quad\quad\quad\quad\cdot\E\big[\veps_{j_1\,i_1}\veps_{j_1\,i_2}
  \cdots \veps_{j_{2k}\,s+i_{2k}}\veps_{j_{2k}\,s+i_1}\big]~.
\end{align}

\smallskip
\noindent Then we assert a lemma stating that $|S(t,s)|\rightarrow 0$ except for one particular term.
\begin{lemma}\label{lemma1}
  $|S(t,s)|\rightarrow 0$ as $p\rightarrow \infty$ unless $t=k$ and $s=k+1$.
\end{lemma}

\noindent Suppose Lemma \ref{lemma1} holds true for a moment, 
then according to \eqref{emka} and \eqref{sts}, we have
\begin{align}\label{f}
  \E m_k(A)=S(k,k+1)+o(1)=\E [\cdot]\cdot \#\{M(A(k,k+1))\}+o(1)~,
\end{align}
where $\E [\cdot]$ refers to the expectation part in \eqref{sts} and $\#\{M(A(k,k+1))\}$ refers to the number of isomorphism class that have $k$ distinct $I$-vertexes and $k+1$ distinct $J$-vertexes.

\noindent First, we show the expectation part $\E [\cdot]$ equals $1$ when $t=k$ and $s=k+1$.
Let $v_m$ denote the number of edges in $M(A(t,s))$ whose degree is $m$. Then we have the total number of  edges having the following relationship:
\begin{align}\label{relation1}
  v_1+2v_2+\cdots +4kv_{4k}=4k~.
\end{align}
Since  we have $\E \veps_{ij}=0$ in \eqref{jus1}, all the multiplicities of the edges in the graph $M(A(t,s))$ should be at least two, that is $v_1=0$.
On the other hand, $M(A(t,s))$ is a connected graph with $t+s$ vertexes and $v_1+\cdots +v_{4k}$ ($=v_2+\cdots +v_{4k}$) edges, we have when $t=k$ and $s=k+1$:
\begin{align}\label{relation2}
  2k+1=t+s&\leq v_1+ \cdots +v_{4k}+1=v_2+ \cdots +v_{4k}+1\nonumber\\
  &\leq \frac 12(2v_2+3v_3+\cdots+4kv_{4k})+1=2k+1~,
\end{align}
where the last equality is due to \eqref{relation1} with $v_1=0$. Then we have all the inequalities in \eqref{relation2} become equalities, that is,
\begin{align*}
  v_2+ \cdots +v_{4k}+1= \frac 12(2v_2+3v_3+\cdots+4kv_{4k})+1=2k+1~,
\end{align*}
which leads to the fact that
\begin{align}\label{edge}
  v_3=v_4=\cdots=v_{4k}=0, \quad v_2=2k~.
\end{align}
This means that all the edges in the graph $M(A(k,k+1))$ is repeated exactly twice, so the part of expectation
\begin{align}\label{expc}
  \E\big[\veps_{j_1\,i_1}\veps_{j_1\,i_2}
  \cdots \veps_{j_{2k}\,s+i_{2k}}\veps_{j_{2k}\,s+i_1}\big]=\left(\E \veps^2_{ji}\right)^{2k}=1~.
\end{align}

\noindent Second, the number of isomorphism class in $M(A(t,s))$ (with each edge repeated at least twice in the original graph $Q(i,j)$) is given by the notation $f_{t-1}(k)$ in \citet{WY14}, where
\begin{align*}
  f_{t-1}(k)=\frac 1k\begin{pmatrix}
    2k \\
    t-1 \\
  \end{pmatrix}\begin{pmatrix}
    k \\
    t \\
  \end{pmatrix}~.
\end{align*}
Therefore, in this special case when $t=k$ and $s=k+1$, we have
\begin{align}\label{makk1}
  \#\{M(A(k,k+1))\}=f_{k-1}(k)=\frac 1k\begin{pmatrix}
    2k \\
    k-1 \\
  \end{pmatrix}~.
\end{align}

\noindent Finally, combine \eqref{f}, \eqref{expc} and \eqref{makk1}, we have
\begin{align*}
  \E m_k(A)=\frac 1k\begin{pmatrix}
    2k \\
    k-1 \\
  \end{pmatrix}+o(1)~.
\end{align*}
Assertion (I) is then proved.
\medskip

\noindent It remains to prove  Lemma \ref{lemma1}.
\begin{proof} (of Lemma \ref{lemma1}) Denote $b_l$ as the degree that associated to the $I$-vertex $i_l$ $(1\leq l\leq t)$ in $M(A(t,s))$, then we have $b_1+\cdots+b_t=4k$, which is the total number of edges. On the other hand, since each edge in $M(A(t,s))$ is repeated at least twice (otherwise, there exist at least one single edge, so the  expectation will be zero), we have each degree $b_l$ at least four (we glue the original $I$-vertexes $i_l$ and $i_l+s$ in $M(A(t,s))$).
  Therefore, we have
  \begin{align*}
    4k=b_1+\cdots+b_t\geq 4t~,
  \end{align*}
  which is $t\leq k$.

  \noindent Now, consider the following two cases separately.
  \medskip

  \noindent {\bf Case 1: $s>k+1$. }

  Recall the definition of $v_m$ in \eqref{relation1}, which satisfies that
  \begin{align*}
    v_1+2v_2+\cdots +4kv_{4k}=2v_2+\cdots +4kv_{4k}=4k~
  \end{align*}
  and
  \begin{align*}
    t+s\leq v_1+\cdots+v_{4k}+1=v_2+\cdots+v_{4k}+1~.
  \end{align*}
  We can bound the expectation part as follows:
  \begin{align}\label{est1}
    &\quad\Big|\E[\veps_{j_1\,i_1}\veps_{j_1\,i_2}\veps_{j_2\,s+i_2}\veps_{j_2\,s+i_3}
    \cdots \veps_{j_{2k}\,s+i_{2k}}\veps_{j_{2k}\,s+i_1}]\Big|\nonumber\\
    &\leq \big|\E \veps^2_{ji}\big|^{v_2}\cdots \big|\E \veps^{4k}_{ji}\big|^{v_{4k}}
    \leq\left(\eta T^{1/4}\right)^{v_3+2v_4+\cdots+(4k-2)v_{4k}}\nonumber\\
    &=\left(\eta T^{1/4}\right)^{3v_3+4v_4+\cdots+4kv_{4k}-2(v_3+v_4+\cdots+v_{4k})}\nonumber\\
    &=\left(\eta T^{1/4}\right)^{4k-2(v_2+v_3+\cdots+v_{4k})}\leq \left(\eta T^{1/4}\right)^{4k-2(t+s-1)}~.
  \end{align}

  Then we have according to \eqref{sts} that
  \begin{align}\label{appe1}
    |S(t,s)|&\leq \frac{1}{p^{k+1}T^k}T^tp^s\left(\eta T^{1/4}\right)^{4k-2(t+s-1)}\#\{M(A(t,s))\}\nonumber\\
    &=\frac{p^{s-k-1}}{T^{\frac 12(s-t-1)}}\eta^{4k-2(t+s-1)}\#\{M(A(t,s))\}\nonumber\\
    &=O\left(\frac{p^{s-k-1}}{T^{\frac 12(s-t-1)}}\eta^{4k-2(t+s-1)}\right)~,
  \end{align}
  where the last equality is due to the fact that $\#\{M(A(t,s))\}$ is a function of $k$ ($k$ is fixed), which could be bounded by a large enough constant.

  Since $s>k+1$ and $t+s-1\leq 2k$, then
  \begin{align*}
    s-k-1-\frac s2+\frac t2+\frac12=\frac s2-k+\frac t2-\frac 12=\frac 12(s+t-2k-1)\leq 0~,
  \end{align*}
  which is
  \begin{align*}
    0<s-k-1\leq \frac 12(s-t-1)~.
  \end{align*}
  So, \eqref{appe1} reduces to
  \begin{align}\label{term1}
    |S(t,k)|\leq O\left(\left(\frac pT\right)^{\frac{s-k-1}{2}}\eta^{4k-2(t+s-1)}\right)\rightarrow 0~,
  \end{align}
  which is due to the fact that $s-k-1>0$ and $p/T\rightarrow 0$.
  \medskip

  \noindent {\bf Case 2}: $s\leq k+1$, but not $t=k$ and $s=k+1$.

  For the same reason as before, we have $t$ distinct $I$-vertexes, each degree is at least four, so we have
  another estimation for the expectation part:
  \begin{align}\label{ex2}
    \bigg|\E[\veps_{j_1\,i_1}\veps_{j_1\,i_2}\veps_{j_2\,s+i_2}\veps_{j_2\,s+i_3}
    \cdots \veps_{j_{2k}\,s+i_{2k}}\veps_{j_{2k}\,s+i_1}]\bigg|\leq \left(\eta T^{1/4}\right)^{4k-4t}~.
  \end{align}
  Therefore,
  \begin{align}
    |S(t,s)|&\leq \frac{1}{p^{k+1}T^k}T^tp^s\left(\eta T^{1/4}\right)^{4k-4t}\#\{M(A(t,s))\}\nonumber\\
    &=O\left(\frac{\eta^{4k-4t}}{p^{k+1-s}}\right)~,
  \end{align}
  which is also due to the fact that $\#\{M(A(t,s))\}=O(1)$.

  Case 2 contains  three situations:
  \begin{align}\label{term2}
    &(1). ~t=k ~\text{and}~ s<k+1: \displaystyle|S(t,s)|\leq O\left(\frac{1}{p^{k+1-s}}\right)\rightarrow 0~;\nonumber\\
    &(2). ~t<k ~\text{and}~ s=k+1: \displaystyle|S(t,s)|\leq O\left(\eta^{4k-4t}\right)\rightarrow 0~;\nonumber\\
    &(3). ~t<k ~\text{and}~ s<k+1: \displaystyle|S(t,s)|\leq O\left(\frac{\eta^{4k-4t}}{p^{k+1-s}}\right)\rightarrow 0~.
  \end{align}

  Combine \eqref{term1} and \eqref{term2}, we have
  $|S(t,k)|\rightarrow 0$ as $p \rightarrow \infty$ unless
  \[
  \left\{
    \begin{array}{l}
      t=k\\
      s=k+1~.
    \end{array}
  \right.
  \]
\end{proof}

\subsection{Proof of Assertion (II)}\label{ssec3}

Recall
\begin{align}\label{varmk}
  &\quad~\Var (m_k(A))\nonumber\\
  &=\frac{1}{p^{2k+2}T^{2k}}\sum_{\bf{i_1},\bf{j_1},\bf{i_2},\bf{j_2}}\left[\E\left(\veps_{Q(\bf{i_1},\bf{j_1})}\veps_{Q(\bf{i_2},\bf{j_2})}\right)-\E\left(\veps_{Q(\bf{i_1},\bf{j_1})}\right)\E\left(\veps_{Q(\bf{i_2},\bf{j_2})}\right)\right]~.
\end{align}
If $Q(\bf{i_1},\bf{j_1})$ has no edges coincident with edges of $Q(\bf{i_2},\bf{j_2})$, then
\[
\E\left(\veps_{Q(\bf{i_1},\bf{j_1})}\veps_{Q(\bf{i_2},\bf{j_2})}\right)-\E\left(\veps_{Q(\bf{i_1},\bf{j_1})}\right)\E\left(\veps_{Q(\bf{i_2},\bf{j_2})}\right)=0
\]
by independence between $\veps_{Q(\bf{i_1},\bf{j_1})}$ and $\veps_{Q(\bf{i_2},\bf{j_2})}$. If $Q=Q({\bf{i_1},\bf{j_1}})\bigcup Q({\bf{i_2},\bf{j_2}})$ has an overall single edge, then
\[
\E\left(\veps_{Q(\bf{i_1},\bf{j_1})}\veps_{Q(\bf{i_2},\bf{j_2})}\right)=\E\left(\veps_{Q(\bf{i_1},\bf{j_1})}\right)\E\left(\veps_{Q(\bf{i_2},\bf{j_2})}\right)=0~,
\]
so in the above two cases, we have
$\Var (m_k(A))=0$.

\noindent Now, suppose $Q=Q({\bf{i_1},\bf{j_1}})\bigcup Q({\bf{i_2},\bf{j_2}})$ has no single edge,  $Q({\bf{i_1},\bf{j_1}})$ and  $Q({\bf{i_2},\bf{j_2}})$ have common edges. Let the number of vertexes of $Q({\bf{i_1},\bf{j_1}})$, $Q({\bf{i_2},\bf{j_2}})$, $Q=Q({\bf{i_1},\bf{j_1}})\bigcup Q({\bf{i_2},\bf{j_2}})$ on the $I$-line be $t_1$, $t_2$, $t$, respectively; and the number of vertexes on the $J$-line be $s_1$, $s_2$, $s$, respectively. Since $Q({\bf{i_1},\bf{j_1}})$ and  $Q({\bf{i_2},\bf{j_2}})$ have common edges, we must have $t\leq t_1+t_2-1$, $s\leq s_1+s_2-1$.

\smallskip
\noindent Similar to \eqref{est1} and \eqref{ex2}, we have two bounds for $\left|\E\left(\veps_{Q(\bf{i_1},\bf{j_1})}\veps_{Q(\bf{i_2},\bf{j_2})}\right)\right|$:
\begin{align}\label{eq1}
  \left|\E\left(\veps_{Q(\bf{i_1},\bf{j_1})}\veps_{Q(\bf{i_2},\bf{j_2})}\right)\right|\leq \left(\eta T^{1/4}\right)^{8k-2(t+s-1)}~,
\end{align}
or
\begin{align}\label{eq2}
  \left|\E\left(\veps_{Q(\bf{i_1},\bf{j_1})}\veps_{Q(\bf{i_2},\bf{j_2})}\right)\right|\leq \left(\eta T^{1/4}\right)^{8k-4t}~.
\end{align}
For the same reason, we have also
\begin{align}\label{eq3}
  \left|\E\veps_{Q(\bf{i_1},\bf{j_1})}\E\veps_{Q(\bf{i_2},\bf{j_2})}\right|&\leq \left(\eta T^{1/4}\right)^{4k-2(t_1+s_1-1)+4k-2(t_2+s_2-1)}\nonumber\\
  &<\left(\eta T^{1/4}\right)^{8k-2(t+s-1)}~,
\end{align}
or
\begin{align}\label{eq4}
  \left|\E\left(\veps_{Q(\bf{i_1},\bf{j_1})}\veps_{Q(\bf{i_2},\bf{j_2})}\right)\right|\leq \left(\eta T^{1/4}\right)^{4k-4t_1+4k-4t_2}<\left(\eta T^{1/4}\right)^{8k-4t}~,
\end{align}
where the last inequalities in \eqref{eq3} and \eqref{eq4} are due to the fact that $t\leq t_1+t_2-1$, $s\leq s_1+s_2-1$.

\smallskip
\noindent Since
\begin{align}\label{vv}
  &\quad~\Var (m_k(A))\nonumber\\
  &=\frac{1}{p^{2k+2}T^{2k}}\sum_{t,s}\sum_{M(A(t,s))}\left[\E\left(\veps_{Q(\bf{i_1},\bf{j_1})}\veps_{Q(\bf{i_2},\bf{j_2})}\right)-\E\left(\veps_{Q(\bf{i_1},\bf{j_1})}\right)\E\left(\veps_{Q(\bf{i_2},\bf{j_2})}\right)\right]\nonumber\\
  &:=\sum_{t,s}\widetilde{S}(t,s)~.
\end{align}

\noindent Using \eqref{eq1}, \eqref{eq2}, \eqref{eq3} and \eqref{eq4}, we can bound the value of $|\widetilde{S}(t,s)|$ as follows:
\begin{align}\label{ss1}
  |\widetilde{S}(t,s)|&\leq O\left(\frac{T^tp^s}{p^{2k+2}T^{2k}}\left(\eta T^{1/4}\right)^{8k-2(t+s-1)}\right)\nonumber\\
  &=O\left(\frac{p^{s-2k-2}}{T^{s/2-t/2-1/2}}\right)~,
\end{align}
or \begin{align}\label{ss2}
  |\widetilde{S}(t,s)|&\leq O\left(\frac{T^tp^s}{p^{2k+2}T^{2k}}\left(\eta T^{1/4}\right)^{8k-4t}\right)\nonumber\\
  &=O\left(p^{s-2k-2}\right)~.
\end{align}

\noindent Clearly,
\begin{align*}
  t_1+s_1\leq 2k+1~, \quad t_2+s_2\leq 2k+1~;
\end{align*}
we have thus
\begin{align*}
  t+s\leq t_1+t_2-1+s_1+s_2-1 \leq 4k~.
\end{align*}

\noindent First, consider the case that $s>t+1$ where we use the bound in \eqref{ss1}.
Since
\begin{align*}
  s-2k-2-s/2+t/2+1/2=s/2+t/2-2k-3/2\leq -3/2~,
\end{align*}
which leads to
\begin{align*}
  s-2k-2\leq -3/2+s/2-t/2-1/2~.
\end{align*}
Combine with \eqref{ss1}, we have
\begin{align}\label{vv1}
  |\widetilde{S}(t,s)|&\leq O\left(p^{-3/2}\left(\frac pT\right)^{\frac 12(s-t-1)}\right)\leq O\left(p^{-3/2}\right)~.
\end{align}

\noindent Second, we use the bound in \eqref{ss2} for the case $s\leq t+1$.
Recall that $t+s\leq 4k$, we have
\begin{align*}
  s-1+s\leq t+s\leq 4k~,
\end{align*}
which is
\begin{align*}
  2s-1\leq 4k~.
\end{align*}
Then, from \eqref{ss2},
\begin{align}\label{vv2}
  |\widetilde{S}(t,s)|&\leq O\left(p^{s-2k-2}\right)\leq O\left(p^{\frac{4k+1}{2}-2k-2}\right)=O\left(p^{-3/2}\right)~.
\end{align}
Combine \eqref{vv}, \eqref{vv1} and \eqref{vv2}, we have
\begin{align*}
  \big|\Var (m_k(A))\big|\leq C(k)p^{-3/2}~,
\end{align*}
which is summable with respect to $p$.
Assertion (II) is then proved.

\section{Convergence of the largest eigenvalue of $A$}\label{largest}

In this section, we  aim to show that the largest eigenvalue of $A$
tends to 4 almost surely, which is the right edge of  its LSD.

\begin{theorem}\label{cb}
  Under the same conditions as in Theorem \ref{singularlsd}, with
  $\sup_{it}\E (\veps_{it}^4)<\infty$ in \eqref{four} replaced by
  $\sup_{it}\E (|\veps_{it}|^{4+\nu})<\infty$ for some $\nu>0$,
  the largest eigenvalue of $A$ converges to $4$ almost surely.
\end{theorem}

Recall that in the proof of Theorem \ref{th1}, a main step is Lemma \ref{lemma1}, which says that $|S(t,s)|\rightarrow 0$ except for one term, which is when $t=k$ and $s=k+1$. One thing to mention here is that in order to prove this lemma, $k$ is assumed to be fixed. Then the number of isomorphism class in $M(A(t,s))$ is a function of $k$, thus can be bounded by a large enough constant. So actually, we do not need to know the value of $\#\{M(A(t,s))\}$ exactly. While in the case of deriving the convergence of the largest eigenvalue, $k$ should grow to infinity, so we can not trivially guarantee  that the number of isomorphism class in $M(A(t,s))$ is still of constant order. Therefore, the main task in this section is to bound this value, making $|S(t,s)|$ ($t\neq k$ or $s\neq k+1$) still a smaller order compared with the main term $|S(k,k+1)|$ when $k \rightarrow \infty$.

\begin{proposition}\label{largek}
  Let the conditions in Theorem \ref{singularlsd} hold, with
  $\sup_{it}\E (\veps_{it}^4)<\infty$ in \eqref{four} replaced by
  $\sup_{it}\E (|\veps_{it}|^{4+\nu})<\infty$ for some $\nu>0$,
  and $k=k(p,T)$ is an integer that tends to infinity and  satisfies the following conditions:
  \begin{align}
    \left\{
      \begin{array}{l}
        k/\log p\rightarrow \infty,\\
        kp/T \rightarrow 0,\\
        k/p \rightarrow. 0\\
      \end{array}
    \right.
  \end{align}
  Then we have
  \begin{align*}
    \E (m_k(A_T))&=\frac 1k \left(\begin{array}{c}
        2k\\
        k-1
      \end{array}\right)\cdot(1+o_k(1))~.
  \end{align*}
\end{proposition}
\medskip

\noindent Now suppose the above Proposition \ref{largek} holds true. We first show it will lead to Theorem \ref{cb}.
\begin{proof}(of Theorem \ref{cb})
  Using Proposition \ref{largek}, we have the estimation that
  \begin{align}\label{estimation}
    \E (m_k(A))=\frac 1k \left(\begin{array}{c}
        2k\\
        k-1
      \end{array}\right)\cdot(1+o_k(1))~,
  \end{align}
  then for any $\Delta>0$, we have
  \begin{align}\label{lb}
    &\quad ~P(l_1>4+\Delta)\leq P(\tr A^k\geq (4+\Delta)^k)\leq\frac{\E \tr A^k}{(4+\Delta)^k}
    =\frac{p\cdot\E(m_k(A))}{(4+\Delta)^k}\nonumber\\
    &\leq \frac{p}{(4+\Delta)^k}\cdot\frac 1k \left(\begin{array}{c}
        2k\\
        k-1
      \end{array}\right)\cdot(1+o_k(1))\leq \left(\frac{4p^{1/k}}{4+\Delta}\right)^k\cdot (1+o_k(1))~.
  \end{align}
  The right hand side tends to $\left(\frac{4}{4+\Delta}\right)^k$ since  $k/\log p\rightarrow \infty$ (so $p^{1/k}\rightarrow 1$). Once we fix this $\Delta>0$,  \eqref{lb} is summable.

  \noindent The upper bound for $l_1$ is trivial due to  our Theorem \ref{th1}.
\end{proof}

\noindent Now it remains to prove our Proposition \ref{largek}.
\begin{proof}(of Proposition \ref{largek})
  After truncation, centralisation and rescaling, we may assume that the $\veps_{it}$'s satisfy the condition that
  \begin{align}\label{cone}
    \ E (\veps_{it})=0,~\Var (\veps_{it})=1,~|\veps_{it}|\leq \delta T^{1/2}~,
  \end{align}
  where $\delta$ is chosen such that
  \begin{align}
    \left\{
      \begin{array}{l}
        \delta \rightarrow 0\\
        \delta T^{1/2-\epsilon}\rightarrow 0\\
        \delta T^{1/2}\rightarrow \infty\\
        \delta^2k\sqrt T\rightarrow 0\\
        \frac{kp}{\delta^2 T}\rightarrow \infty~.
      \end{array}
    \right.
  \end{align}
  More detailed justifications of \eqref{cone} are provided in  Appendix \ref{lartr}.

  From the proof of Theorem \ref{th1}, we have
  \begin{align*}
    \E m_k(A)=\sum_{t,s}S(t,s)=S(k,k+1)+o(1)=\frac 1k \begin{pmatrix}
      2k \\
      k-1 \\
    \end{pmatrix}+o(1)~,
  \end{align*}
  where $S(k,k+1)$ is the main term that contributes to $\E m_k(A)$, while all other terms can be neglect.
  Therefore, it remains to prove that when $k \rightarrow \infty$, we still have
  \begin{align*}
    \sum_{t\neq k ~\text{or}~ s\neq k+1}S(t,s)=\frac 1k \begin{pmatrix}
      2k \\
      k-1 \\
    \end{pmatrix}\cdot o_k(1) ~.
  \end{align*}
  We also consider two cases:

  \begin{center}
    $
    \begin{array}{l}
      \text{Case 1}: s>k+1 \\[2mm]
      \text{Case 2}: s\leq k+1, ~\text{but not}~ t=k ~\text{and}~ s=k+1.
    \end{array}
    $
  \end{center}

  \noindent Similar to \eqref{est1} and \eqref{ex2}, we have two bounds for the expectation part:
  \begin{align}\label{new1}
    \Big|\E[\veps_{j_1\,i_1}\veps_{j_1\,i_2}\veps_{j_2\,s+i_2}\veps_{j_2\,s+i_3}
    \cdots \veps_{j_{2k}\,s+i_{2k}}\veps_{j_{2k}\,s+i_1}]\Big|
    \leq \left(\delta T^{1/2}\right)^{4k-2(t+s-1)}~
  \end{align}
  or
  \begin{align}\label{new2}
    \bigg|\E[\veps_{j_1\,i_1}\veps_{j_1\,i_2}\veps_{j_2\,s+i_2}\veps_{j_2\,s+i_3}
    \cdots \veps_{j_{2k}\,s+i_{2k}}\veps_{j_{2k}\,s+i_1}]\bigg|\leq \left(\delta T^{1/2}\right)^{4k-4t}~.
  \end{align}

  \noindent Consider  $t=1$ first. From \citet{wang}, the number of isomorphism class $\#\{M(A(1,s))\}$ is bounded by
  \[\begin{pmatrix}
    2k \\
    2k-s \\
  \end{pmatrix}~,
  \]
  and combine this with \eqref{sts} and \eqref{new1}, we have
  \begin{align}
    |S(1,s)|\leq \frac{1}{p^{k+1}T^k}Tp^s\left(\delta T^{1/2}\right)^{4k-2s}\left(\begin{array}{c}
        2k\\
        2k-s
      \end{array}\right)~.
  \end{align}
  Then,
  \begin{align}\label{s1}
    \big|\sum_{s}S(1,s)\big|&\leq \sum_{s=1}^{2k}\frac{1}{p^{k+1}T^k}Tp^s\left(\delta T^{1/2}\right)^{4k-2s}\left(\begin{array}{c}
        2k\\
        2k-s
      \end{array}\right)\nonumber\\
    &= \sum_{s=1}^{2k}\frac{1}{p^{k+1}T^k}Tp^s\left(\delta T^{1/2}\right)^{4k-2s}\left(\begin{array}{c}
        2k\\
        s
      \end{array}\right)~.
  \end{align}
  The right hand side of \eqref{s1} can be bounded as
  \begin{align*}
    \sum_{s=1}^{2k}\frac{T}{p^{k+1}T^k}\left(\delta T^{1/2}\right)^{4k}\left(\frac{2kp}{\delta^2  T}\right)^s~,
  \end{align*}
  which is dominated by the term when $s=2k$ since $\frac{kp}{\delta^2  T}\rightarrow \infty$. Then \eqref{s1} reduces to
  \begin{align}
    \frac{1}{p^{k+1}T^k}Tp^{2k}=\left(\frac p T\right)^{k-1}\rightarrow 0~.
  \end{align}

  \noindent Next, we  consider Case 1 and Case 2 (when $t>1$) separately. According to \citet{wang}, the number of isomorphism class in $M(A(s,t))$ ($t>1$) is bounded by
  \begin{align}\label{iso}
    f_{t-1}(k)\left(\begin{array}{c}
        2k-t\\
        s-1
      \end{array}\right)~,
  \end{align}
  where
  \begin{align*}
    f_{t-1}(k)=\frac 1k\begin{pmatrix}
      2k \\
      t-1 \\
    \end{pmatrix}\begin{pmatrix}
      k \\
      t \\
    \end{pmatrix}~.
  \end{align*}
  Case 1 ($s>k+1$ and $t>1$):
  The part of expectation can be bounded by \eqref{new1}, and
  combining this with \eqref{sts} and \eqref{iso}, we have
  \begin{align}
    |S(t,s)|&\leq\frac{1}{p^{k+1}T^k}\sum_{M(A(t,s))}p^sT^t
    \cdot\Big|\E[\veps_{j_1\,i_1}\veps_{j_1\,i_2}
    \cdots \veps_{j_{2k}\,s+i_{2k}}\veps_{j_{2k}\,s+i_1}]\Big|\nonumber\\
    &\leq \frac{p^sT^t}{p^{k+1}T^k}\left(\delta T^{1/2}\right)^{4k-2(t+s-1)}\cdot  f_{t-1}(k)\left(\begin{array}{c}
        2k-t\\
        s-1
      \end{array}\right)~.
  \end{align}
  Since  $s\geq k+2$, $t\geq 2$, and a trivial relationship  that $t+s-1\leq 2k$, we have
  \begin{align}\label{em1}
    \Big|\sum_{t,s}S(t,s)\Big|&\leq \sum_{t=2}^{k-1}\sum_{s=k+2}^{2k+1-t}\frac{p^sT^t}{p^{k+1}T^k}\left(\delta T^{1/2}\right)^{4k-2(t+s-1)}\cdot f_{t-1}(k)\left(\begin{array}{c}
        2k-t\\
        s-1
      \end{array}\right)~.
  \end{align}
  The summation over $s$ in \eqref{em1} can be bounded as follows:
  \begin{align}\label{em2}
    \sum_{s=k+2}^{2k+1-t}\delta^{-2s}T^{-s}p^s\left(\begin{array}{c}
        2k-t\\
        s-1
      \end{array}\right)\leq \sum_{s=k+2}^{2k+1-t}\left(\frac{2kp}{\delta^2 T}\right)^s~,
  \end{align}
  and since $\frac{kp}{\delta^2 T}\rightarrow \infty$, the summation in \eqref{em2} is dominated by the term of $s=2k+1-t$. Therefore, \eqref{em1} reduces to
  \begin{align}\label{em3}
    &\quad \sum_{t=2}^{k-1}\frac{p^{2k+1-t}T^t}{p^{k+1}T^k}\left(\delta T^{1/2}\right)^{4k-2(t+2k+1-t-1)}\cdot f_{t-1}(k)\left(\begin{array}{c}
        2k-t\\
        2k+1-t-1
      \end{array}\right)\nonumber\\
    &=\sum_{t=2}^{k-1} \left(\frac{p}{T}\right)^{k-t}f_{t-1}(k)=\sum_{t=2}^{k-1} \frac 1k \left(\begin{array}{c}
        2k\\
        t-1\end{array}\right)\left(\begin{array}{c}k\\
        t\end{array}\right)\left(\frac{p}{T}\right)^{k-t}~.
  \end{align}

  \noindent For the same reason, the right hand side of \eqref{em3} inside the summation can be bounded by
  \[\frac 1k\left(\frac pT\right)^k\left(\frac{2k^2T}{p}\right)^t~,\]
  and since $Tk^2/p=k^2/(\frac p T) \rightarrow \infty$, the dominating term in \eqref{em3} is when $t=k-1$, which reduces to
  \begin{align}\label{em4}
    \frac 1k \left(\begin{array}{c}
        2k\\
        k-2\end{array}\right)\left(\begin{array}{c}k\\
        k-1\end{array}\right)\left(\frac{p}{T}\right)=\frac{k(k-1)}{k+2}\frac pT\cdot \frac 1k\left(\begin{array}{c}2k\\k-1\end{array}\right)~.
  \end{align}
  Since $kp/T\rightarrow 0$, we have \eqref{em4} equals
  \begin{align}
    \frac 1k\left(\begin{array}{c}2k\\k-1\end{array}\right)\cdot o_k(1)~.
  \end{align}
  Therefore, in this case, we have
  \begin{align}
    \Big|\sum_{t,s}S(t,s)\Big|=\frac 1k\left(\begin{array}{c}2k\\k-1\end{array}\right)\cdot o_k(1)~.
  \end{align}

  \medskip
  \noindent Case 2 ($2\leq t\leq k$ and $s\leq k+1$):
  For the same reason, combining the bound of the expectation part in \eqref{new2} with \eqref{sts} and \eqref{iso}, we have
  \begin{align}
    |S(t,s)|&=\frac{1}{p^{k+1}T^k}\sum_{M(A(t,s))}p^sT^t\cdot\Big|\E[\veps_{j_1\,i_1}\veps_{j_1\,i_2}
    \cdots \veps_{j_{2k}\,s+i_{2k}}\veps_{j_{2k}\,s+i_1}]\Big|\nonumber\\
    &\leq \frac{1}{p^{k+1}T^k}\left(\delta T^{1/2}\right)^{4k-4t}\cdot p^sT^t\cdot f_{t-1}(k)\left(\begin{array}{c}
        2k-t\\
        s-1
      \end{array}\right)~.
  \end{align}
  Therefore, we have
  \begin{align}\label{em6}
    \Big|\sum_{t,s}S(t,s)\Big|\leq\sum_{t=2}^k \sum_{s=1}^{k+1}p^{s-k-1}T^{k-t}\delta^{4k-4t}f_{t-1}(k)\left(
      \begin{array}{c}
        2k-t\\
        s-1
      \end{array}\right)~.
  \end{align}
  We also consider the following three situations:
  \begin{align*}
    &(1). ~t=k ~\text{and}~ s<k+1~,\\
    &(2). ~1<t<k ~\text{and}~ s=k+1~,\\
    &(3). ~1<t<k ~\text{and}~ s<k+1~,
  \end{align*}
  and show that for all the above three situations, we have \eqref{em6} bounded by
  \begin{align*}
    \frac 1k \begin{pmatrix}
      2k \\
      k-1 \\
    \end{pmatrix}\cdot o_k(1)~.
  \end{align*}
  For situation (1), \eqref{em6} reduces to
  \begin{align}\label{t1}
    \sum_{s=1}^{k}p^{s-k-1}f_{k-1}(k)\left(\begin{array}{c}
        k\\
        s-1
      \end{array}\right)
    =\sum_{s=1}^{k}p^{s-k-1}\frac 1k \left(\begin{array}{c}
        2k\\
        k-1
      \end{array}\right)\left(\begin{array}{c}
        k\\
        s-1
      \end{array}\right)~,
  \end{align}
  which can be bounded as
  \begin{align*}
    \sum_{s=1}^k p^{-k-1}\frac 1k \begin{pmatrix}
      2k \\
      k-1 \\
    \end{pmatrix}(kp)^s~.
  \end{align*}
  Therefore, the dominating term is when $s=k$, thus \eqref{t1} reduces to
  \begin{align*}
    \frac 1k \begin{pmatrix}
      2k \\
      k-1 \\
    \end{pmatrix}\cdot \frac kp=\frac 1k \begin{pmatrix}
      2k \\
      k-1 \\
    \end{pmatrix}\cdot o_k(1)~,
  \end{align*}
  which is due to the choice of $k$ that $k/p \rightarrow 0$.
  \smallskip

  \noindent For situation (2), \eqref{em6} reduces to
  \begin{align}\label{en88}
    &~\quad\sum_{t=2}^{k-1} \delta^{4k-4t}T^{k-t}\cdot f_{t-1}(k)\left(\begin{array}{c}
        2k-t\\
        k
      \end{array}\right)\nonumber\\
    &=\sum_{t=2}^{k-1}\delta^{4k-4t}T^{k-t}\cdot\frac 1k \left(\begin{array}{c}
        2k\\
        t-1
      \end{array}\right)\left(\begin{array}{c}
        k\\
        t
      \end{array}\right)\left(\begin{array}{c}
        2k-t\\
        k
      \end{array}\right)~.
  \end{align}

  \noindent Since the right hand side of \eqref{en88} can be bounded by
  \begin{align}
    \sum_{t=2}^{k-1}\delta^{4k}\cdot\frac{(2kT)^k}{k}\left(\frac{2k^2}{\delta^4T}\right)^t~,
  \end{align}
  which is dominated by the term of $t=k-1$ since $\frac{2k^2}{T\delta^4}=\frac{2k^2}{(\delta T^{1/4})^4}\rightarrow \infty$. Therefore, we have \eqref{en88} bounded by
  \begin{align*}
    &~\quad T \delta^{4}\cdot\frac 1k \left(\begin{array}{c}
        2k\\
        k-2
      \end{array}\right)\left(\begin{array}{c}
        k\\
        k-1
      \end{array}\right)\left(\begin{array}{c}
        k+1\\
        k
      \end{array}\right)~\nonumber\\
    &=\delta^4 T\frac{(k-1)k(k+1)}{k+2}\cdot\frac 1k \left(\begin{array}{c}
        2k\\
        k-1
      \end{array}\right)\nonumber\\
    &\rightarrow \frac 1k \left(\begin{array}{c}
        2k\\
        k-1
      \end{array}\right)\cdot o_k(1)~,
  \end{align*}
  which is due to the fact that $\delta^4 Tk^2=(\delta^2 \sqrt T k)^2\rightarrow 0$.

  \smallskip

  \noindent For situation (3), we have \eqref{em6} reduce to
  \begin{align}\label{en7}
    &~\quad\sum_{t=2}^{k-1}\sum_{s=1}^kp^{s-k-1}T^{k-t} \delta^{4k-4t}\cdot f_{t-1}(k)\left(\begin{array}{c}
        2k-t\\
        s-1
      \end{array}\right)\nonumber\\
    &=\sum_{t=2}^{k-1}\sum_{s=1}^kp^{s-k-1}T^{k-t}\delta^{4k-4t}\cdot\frac 1k \left(\begin{array}{c}
        2k\\
        t-1
      \end{array}\right)\left(\begin{array}{c}
        k\\
        t
      \end{array}\right)\left(\begin{array}{c}
        2k-t\\
        s-1
      \end{array}\right)~.
  \end{align}
  The part of summation over $s$ is
  \begin{align*}
    \sum_{s=1}^k p^s\left(\begin{array}{c}
        2k-t\\
        s-1
      \end{array}\right)~,
  \end{align*}
  which could be bounded by
  \begin{align*}
    \sum_{s=1}^k (2kp)^s~,
  \end{align*}
  therefore, the dominating term is when $s=k$. So \eqref{en7} reduces to
  \begin{align}\label{en8}
    \sum_{t=2}^{k-1}p^{-1}\delta^{4k-4t}T^{k-t}\cdot\frac 1k \left(\begin{array}{c}
        2k\\
        t-1
      \end{array}\right)\left(\begin{array}{c}
        k\\
        t
      \end{array}\right)\left(\begin{array}{c}
        2k-t\\
        k-1
      \end{array}\right)~.
  \end{align}
  For the same reason, the right hand side of \eqref{en8} can be bounded by
  \begin{align*}
    \sum_{t=2}^{k-1}p^{-1}\delta^{4k}\cdot\frac 1k\left(\frac{2k^2}{T\delta^{4}}\right)^t(2kT)^k~,
  \end{align*}
  which is dominated by the term of $t=k-1$ since $\frac{k^2}{T\delta^4}=\frac{k^2}{\left(\delta T^{1/4}\right)^4}\rightarrow \infty$. Therefore, \eqref{en8} reduces to
  \begin{align}\label{en9}
    \frac{T}{p}\delta^{4}\cdot\frac 1k \left(\begin{array}{c}
        2k\\
        k-2
      \end{array}\right)\left(\begin{array}{c}
        k\\
        k-1
      \end{array}\right)\left(\begin{array}{c}
        k+1\\
        k-1
      \end{array}\right)
    =O\left(\frac{\delta^4 k^3T}{p}\cdot\frac 1k \left(\begin{array}{c}2k\\
          k-1\end{array}\right)\right)~,
  \end{align}
  and since $\frac{\delta^4 k^3T}{p}=(\delta^2 k\sqrt T)^2\cdot k/p \rightarrow 0$, we have \eqref{en9} equals
  \begin{align*}
    \frac 1k \left(\begin{array}{c}2k\\
        k-1\end{array}\right)\cdot o_k(1)~.
  \end{align*}
  Finally, in all the three situations, we  have
  \begin{align*}
    \Big|\sum_{t,s}S(t,s)\Big|=\frac 1k \left(\begin{array}{c}2k\\
        k-1\end{array}\right)\cdot o_k(1)~.
  \end{align*}

  \noindent The proof of Proposition \ref{largek} is complete.

\end{proof}

\appendix

\section{Justification of truncation, centralisation and rescaling in \eqref{jus1}}\label{j1}
\subsection{Truncation}
Define two $p \times T$ matrices
\begin{align}\label{e1e2}
  &E_1:=(\veps_1 ~ \veps_2\cdots ~ \veps_{T-1} ~ \veps_{T})~,\quad
  E_2:=(\veps_{s+1} ~ \veps_{s+2}\cdots ~ \veps_{s+T-1} ~ \veps_{s+T})~,
\end{align}
then
\begin{align}\label{xtt}
  X_T=\frac 1T \sum_{t=s+1}^{s+T}\veps_t \veps^T_{t-s}=\frac 1T E_2E_1^T~,
\end{align}
and our target matrix
\begin{align}\label{A}
  A=\frac Tp X_TX_T^T=\frac {1}{pT}E_2E_1^TE_1E_2^T~.
\end{align}
Let
\begin{align*}
  \hat{\veps}_{ij}=\veps_{ij}\mathbf{1}_{\{|\veps_{ij}|\leq \eta T^{1/4}\}}~,
\end{align*}
$\hat{X}_T$ and $\hat{A}$ are defined by replacing all the $\veps_{ij}$ with $\hat{\veps}_{ij}$ in \eqref{xtt} and \eqref{A}.

\noindent Using Theorem A.44 in \citet{bai} and the inequality that
\[
\text{rank} (AB-CD)\leq \text{rank}(A-C)+\text{rank}(B-D)~,
\]
we have
\begin{align}\label{tru}
  &\quad\left\|F^A(x)-F^{\hat{A}}(x)\right\|=\left\|F^{\frac T p X_TX^T_T}(x)-F^{\frac T p \hat{X}_T\hat{X}^T_T}(x)\right\|\nonumber\\
  &\leq \frac 1p \text{rank}\left(\sqrt{\frac T p}X_T-\sqrt{\frac T p}\hat{X}_T\right)
  = \frac 1p \text{rank}\left(X_T-\hat{X}_T\right)\nonumber\\
  &=\frac 1p \text{rank} \left(\frac 1T E_2E^T_1-\frac 1T \hat{E}_2\hat{E}^T_1\right)=\frac 1p \text{rank} \left( E_2E^T_1- \hat{E}_2\hat{E}^T_1\right)\nonumber\\
  &\leq \frac 1p \text{rank} \left(E_2-\hat{E}_2\right)+\frac 1p \text{rank} \left(E_1-\hat{E}_1\right)\nonumber\\
  &=\frac 2p \text{rank} \left(E_1-\hat{E}_1\right)\leq \frac 2p \sum_{i=1}^p \sum_{j=1}^T \mathbf{1}_{\{|\veps_{ij}|>\eta T^{1/4}\}}~.
\end{align}

\noindent Since
$\sup_{it}\E (\veps_{it}^4)<\infty$, we have always
\[
\frac{1}{\eta^4 pT}\sum_{i,j}\E \left(|\veps_{ij}|^4I_{(|\veps_{ij}|> \eta T^{1/4})}\right)\longrightarrow 0~\quad \text{as}~ p, T \rightarrow \infty~.
\]

\noindent Consider the expectation and variance of $\frac 1p \sum_{i=1}^p \sum_{j=1}^T \mathbf{1}_{\{|\veps_{ij}|>\eta T^{1/4}\}}$ in \eqref{tru}:
\begin{align*}
  \E\bigg(\frac 2p \sum_{i=1}^p \sum_{j=1}^T \mathbf{1}_{\{|\veps_{ij}|>\eta T^{1/4}\}}\bigg)&\leq
  \frac 2p \sum_{i=1}^p \sum_{j=1}^T \frac{\E \Big(|\veps_{ij}|^4\cdot\mathbf{1}_{\{|\veps_{ij}|>\eta T^{1/4}\}}\Big)}{\eta^4 T}=o(1)~,\\
  \Var\bigg( \frac 2p \sum_{i=1}^p \sum_{j=1}^T \mathbf{1}_{\{|\veps_{ij}|>\eta T^{1/4}\}}\bigg)&\leq
  \frac {4}{p^2} \sum_{i=1}^p \sum_{j=1}^T \frac{\E \Big(|\veps_{ij}|^4\cdot\mathbf{1}_{\{|\veps_{ij}|>\eta T^{1/4}\}}\Big)}{\eta^4 T}=o(\frac 1p)~.
\end{align*}

\noindent Applying Bernstein's inequality, for all small $\veps>0$ and large $p$, we have
\begin{align}\label{pro}
  &P\left(\frac 2p \sum_{i=1}^p \sum_{j=1}^T \mathbf{1}_{\{|\veps_{ij}|>\eta T^{1/4}\}}\geq \veps\right)\leq 2e^{-\frac 12 \veps^2 p}~.
\end{align}

\noindent Finally, combine \eqref{tru}, \eqref{pro} with Borel-Cantelli lemma, we have with probability 1,
\begin{align*}
  \left\|F^A(x)-F^{\hat{A}}(x)\right\|\rightarrow 0~.
\end{align*}

\subsection{Centralisation}
Let
\begin{align*}
  \tilde{\veps}_{ij}=\hat{\veps}_{ij}-\E \hat{\veps}_{ij}~,
\end{align*}
$\tilde{X}_T$ and $\tilde{A}$ are  defined by involving the $\tilde{\veps}_{ij}$'s in  \eqref{xtt} and \eqref{A}.
\smallskip

\noindent Similar to \eqref{tru}, we have
\begin{align*}
  &\quad\left\|F^{\hat{A}}(x)-F^{\tilde{A}}(x)\right\|
  \leq \frac 1p \text{rank}\left(\sqrt{\frac T p}\hat{X}_T-\sqrt{\frac T p}\tilde{X}_T\right)\nonumber\\
  &= \frac 1p \text{rank}\left(\hat{E}_2\hat{E}^T_1- \tilde{E}_2\tilde{E}^T_1\right)
  \leq \frac 2p \text{rank} \left(\hat{E}_1-\tilde{E}_1\right)\nonumber\\
  &=\frac 2p \text{rank}\big(\E (\hat{E}_1)\big)=\frac 2p \rightarrow 0~,\quad \text{as}~p \rightarrow \infty~.
\end{align*}
Therefore, we have
\begin{align*}
  \left\|F^{\hat{A}}(x)-F^{\tilde{A}}(x)\right\|\rightarrow 0~.
\end{align*}

\subsection{Rescaling}
Let
\begin{align*}
  \sigma^2_{ij}=\E \tilde{\veps}^2_{ij}~,\quad \check{\veps}_{ij}:=\tilde{\veps}_{ij}/\sigma_{ij}~,
\end{align*}
then for the same reason as \eqref{tru}, we have
\begin{align*}
  &\quad\left\|F^{\tilde{A}}(x)-F^{\check{A}}(x)\right\|
  \leq \frac 1p \text{rank}\left(\tilde{E}_2-\check{E}_2\right)+\frac 1p \text{rank}\left(\tilde{E}_1-\check{E}_1\right)\nonumber\\
  &= \frac 2p \text{rank} \left(\tilde{E}_1-\check{E}_1\right)\leq
  \frac 2p \max_{i\,j} \left(1-\frac{1}{\sigma_{ij}}\right)\cdot\text{rank} \left(\tilde{E}_1\right)\\
  & \leq \frac 2p \max_{i\,j} \left(1-\frac{1}{\sigma_{ij}}\right)\min\{p,T\}\\
  & = O\left( \max_{i\,j} \left(1-\frac{1}{\sigma_{ij}}\right)\right)~.
\end{align*}
Since
\begin{align*}
  &\sigma^2_{ij}=\E \tilde{\veps}^2_{ij}=\E (\hat{\veps}_{ij}-\E \hat{\veps}_{ij})^2=\Var(\hat{\veps}_{ij})
  = \Var(\veps_{ij}\cdot\mathbf{1}_{\{|\veps_{ij}|\leq \eta T^{1/4}\}})\\
  &\quad ~\rightarrow \Var(\veps_{ij})=1~, ~\text{as}~T \rightarrow \infty~.
\end{align*}

\noindent Therefore, we have
\[
\left\|F^{\tilde{A}}(x)-F^{\check{A}}(x)\right\|\rightarrow 0~.
\]

\section{Justification of truncation, centralisation and rescaling in \eqref{cone}}\label{lartr}
\subsection{Truncation}
$E_1$, $E_2$, $X_T$ and $A$ are defined in \eqref{e1e2}, \eqref{xtt} and \eqref{A}.
Let
\begin{align*}
  \hat{\veps}_{ij}=\veps_{ij}\mathbf{1}_{\{|\veps_{ij}|\leq \delta T^{1/2}\}}~,
\end{align*}
$\hat{X}_T$ and $\hat{A}$ are defined by replacing all the $\veps_{ij}$ with $\hat{\veps}_{ij}$ in  \eqref{xtt} and \eqref{A}.
With the assumption that
$\sup_{it}\E (|\veps_{it}|^{4+\nu})<\infty$, we have always
\begin{align}\label{condi}
  \sup_{it} \frac{\E\Big(|\veps_{it}|^{4+\nu}\mathbf{1}_{\{|\veps_{it}|> \delta T^{1/2}\}}\Big)}{\delta^{4+\nu} }\longrightarrow 0~\quad \text{as}~ p, T \rightarrow \infty~.
\end{align}
Since
\begin{align*}
  A=\frac {1}{pT}E_2E_1^TE_1E_2^T~,
\end{align*}
whose eigenvalues are the same as those of
\begin{align*}
  B:=\frac {1}{pT}E_1^TE_1E_2^TE_2~,
\end{align*}
then we have
\begin{align}\label{j12}
  &\quad\left|\lambda_{\max}(A)-\lambda_{\max}(\hat{A})\right|=\left|\lambda_{\max}(B)-\lambda_{\max}(\hat{B})\right|\nonumber\\
  &=\left|\left\|\frac {1}{pT}E^T_1E_1E^T_2E_2\right\|_{op}-\left\|\frac {1}{pT}\hat{E}^T_1\hat{E}_1\hat{E}^T_2\hat{E}_2\right\|_{op}\right|\nonumber\\
  & \leq\left\|\frac {1}{pT}E^T_1E_1E^T_2E_2-\frac {1}{pT}\hat{E}^T_1\hat{E}_1\hat{E}^T_2\hat{E}_2\right\|_{op}\nonumber\\
  &\leq\left\|\frac {1}{pT}E^T_1E_1E^T_2E_2-\frac {1}{pT}\hat{E}^T_1\hat{E}_1E^T_2E_2\right\|_{op}+\left\|\frac {1}{pT}\hat{E}^T_1\hat{E}_1E^T_2E_2-\frac {1}{pT}\hat{E}^T_1\hat{E}_1\hat{E}^T_2\hat{E}_2\right\|_{op}\nonumber\\
  &=\left\|\frac {1}{pT}\left(E^T_1E_1-\hat{E}^T_1\hat{E}_1\right)E^T_2E_2\right\|_{op}+\left\|\frac {1}{pT}\hat{E}^T_1\hat{E}_1\left(E^T_2E_2-\hat{E}^T_2\hat{E}_2\right)\right\|_{op}\nonumber\\
  &:=J_1+J_2~.
\end{align}

\noindent First, we have
\begin{align}\label{jjj}
  &\quad\left\|E^T_1E_1-\hat{E}^T_1\hat{E}_1\right\|_{op}=\max_{\|x\|=1}x(E^T_1E_1-\hat{E}^T_1\hat{E}_1)x^T\nonumber\\
  &=\max_{\|x\|=1}\left[x(E^T_1E_1-\hat{E}^T_1E_1)x^T+x(\hat{E}^T_1E_1-\hat{E}^T_1\hat{E}_1)x^T\right]\nonumber\\
  &\leq\max_{\|x\|=1}x(E^T_1E_1-\hat{E}^T_1E_1)x^T+\max_{\|x\|=1}x(\hat{E}^T_1E_1-\hat{E}^T_1\hat{E}_1)x^T\nonumber\\
  &:=J_{11}+J_{12}~,
\end{align}
where
\begin{align}\label{j11}
  J_{11}&=\max_{\|x\|=1}x(E^T_1E_1-\hat{E}^T_1E_1)x^T=\max_{\|x\|=1}\sum_{i,j}x_ix_j(E^T_1E_1-\hat{E}^T_1E_1)(i,j)\nonumber\\
  &=\max_{\|x\|=1}\sum_{i,j}x_ix_j\sum_{k=1}^p\left(\veps_{ki}-\hat{\veps}_{ki}\right)\veps_{kj}\nonumber\\
  & \leq \max_{\|x\|=1}\sum_{k=1}^p \bigg[\Big(\sum_i x^2_i\Big)^{1/2}\Big(\sum_i\big(\veps_{ki}-\hat{\veps}_{ki}\big)^2\Big)^{1/2}\cdot
  \Big(\sum_j x^2_j\Big)^{1/2}\Big(\sum_j\veps^2_{kj}\Big)^{1/2} \bigg]\nonumber\\
  &= \sum_{k=1}^p \bigg[\Big(\sum_i\big(\veps_{ki}-\hat{\veps}_{ki}\big)^2\Big)^{1/2}\cdot
  \Big(\sum_j\veps^2_{kj}\Big)^{1/2} \bigg]\nonumber\\
  &\leq \bigg(\sum_{k=1}^p \sum_{i=1}^T\big(\veps_{ki}-\hat{\veps}_{ki}\big)^2\bigg)^{1/2}\cdot
  \bigg(\sum_{k=1}^p\sum_{j=1}^T\veps^2_{kj}\bigg)^{1/2}\nonumber\\
  &=O\left(\sqrt{pT}\cdot\bigg(\sum_{k=1}^p \sum_{i=1}^T\veps^2_{ki}\cdot\mathbf{1}_{\{|\veps_{ki}|>\delta T^{1/2}\}}\bigg)^{1/2}\right)\nonumber\\
  &\leq O\left(\bigg((pT)^2\cdot \sup_{k\,i}\E \Big(\veps^2_{ki}\cdot \mathbf{1}_{\{|\veps_{ki}|>\delta T^{1/2}\}}\Big)\bigg)^{1/2}\right)\nonumber\\
  &\leq O\left(\bigg(\frac{(pT)^2}{(\delta T^{1/2})^{2+\nu}}\cdot \sup_{k\,i}\E \Big(|\veps_{ki}|^{4+\nu}\cdot \mathbf{1}_{\{|\veps_{ki}|>\delta T^{1/2}\}}\Big)\bigg)^{1/2}\right)\nonumber\\
  &=o\Big(\delta p T^{1/2-\nu/4}\Big)~,
\end{align}
where the last inequality is due to \eqref{condi}.

\noindent For the same reason, $J_{12}$ is also of the same order as \eqref{j11}. Therefore we have
\begin{align}\label{eqa1}
  \quad\left\|E^T_1E_1-\hat{E}^T_1\hat{E}_1\right\|_{op}\leq o\Big(\delta p T^{1/2-\nu/4}\Big)~.
\end{align}
Then recall the definition of $J_1$ in \eqref{j12}, where
\begin{align}\label{j}
  J_1&=\left\|\frac {1}{pT}\left(E^T_1E_1-\hat{E}^T_1\hat{E}_1\right)E^T_2E_2\right\|_{op}\leq \frac {1}{p}\left\|E^T_1E_1-\hat{E}^T_1\hat{E}_1\right\|_{op}\cdot \frac 1T\left\|E^T_2E_2\right\|_{op}\nonumber\\
  &\leq o\Big(\delta  T^{1/2-\nu/4}\Big)\rightarrow 0~, \quad \text{as}~p, T \rightarrow \infty~,
\end{align}
where the last inequality in \eqref{j} is due to \eqref{eqa1} and the fact that
$\frac 1T\left\|E^T_2E_2\right\|_{op}$ is the largest eigenvalue of the sample covariance matrix $\frac 1T E_2E^T_2$, which is of constant order.

\noindent For the same reason, we also have $J_2$ the same order as $J_1$, which also tends to zero. Finally, according to \eqref{j12} we have
\begin{align*}
  \left|\lambda_{\max}(A)-\lambda_{\max}(\hat{A})\right|\rightarrow 0~.
\end{align*}

\subsection{Centralisation and Rescaling}
Let
\begin{align*}
  \sigma^2_{it}=\Var \hat{\veps}_{it}~,\quad \tilde{\veps}_{it}=\frac{\hat{\veps}_{it}-\E \hat{\veps}_{it}}{\sigma_{it}}~,
\end{align*}
$\tilde{X}_T$ and $\tilde{A}$ are defined by replacing all the $\veps_{ij}$ with $\tilde{\veps}_{ij}$ in \eqref{xtt} and \eqref{A}. In this subsection, we will show
\begin{align*}
  \left|\lambda_{\max}(\hat{A})-\lambda_{\max}(\tilde{A})\right|\rightarrow 0~,
\end{align*}
which is equivalent to showing
\begin{align*}
  \left|\lambda_{\max}(\hat{B})-\lambda_{\max}(\tilde{B})\right|\rightarrow 0~.
\end{align*}
First, since
\begin{align}\label{si}
  &~\quad\sup_{ki}\left|1-\sigma^2_{ki}\right|
  =\sup_{ki}\left|\E \veps^2_{ki}-\E \left(\veps_{ki}\cdot \mathbf{1}_{\{|\veps_{ki}|\leq \delta T^{1/2}\}}-\E \Big( \veps_{ki}\cdot \mathbf{1}_{\{|\veps_{ki}|\leq \delta T^{1/2}\}}\Big)\right)^2\right|\nonumber\\
  &=\sup_{ki}\left|\E\Big( \veps^2_{ki}\cdot \mathbf{1}_{\{|\veps_{ki}|> \delta T^{1/2}\}}\Big)+\Big(\E \big(\veps_{ki}\cdot \mathbf{1}_{\{|\veps_{ki}|> \delta T^{1/2}\}}\big)\Big)^2\right|\nonumber\\
  &\leq2\cdot\sup_{ki}\left|\E\Big( \veps^2_{ki}\cdot \mathbf{1}_{\{|\veps_{ki}|> \delta T^{1/2}\}}\Big)\right|
  \leq\frac{2\cdot\sup_{ki}\E\Big( |\veps_{ki}|^{4+\nu}\cdot \mathbf{1}_{\{|\veps_{ki}|> \delta T^{1/2}\}}\Big)}{\big(\delta T^{1/2}\big)^{2+\nu}}\nonumber\\
  &=o\left(\frac{\delta^2}{T^{\frac{2+\nu}{2}}}\right)~,
\end{align}
where the last equality is due to \eqref{condi}.
Finally, we have:
\begin{align}\label{bound1}
  \sup_{ki}\left|1-\frac{1}{\sigma_{ki}}\right|&=\sup_{ki}\left|\frac{\sigma_{ki}-1}{\sigma_{ki}}\right|
  =\sup_{ki}\left|\frac{\sigma^2_{ki}-1}{\sigma_{ki}(\sigma_{ki}+1)}\right|=O\left(\sup_{ki}\left|\sigma^2_{ki}-1\right|\right)\nonumber\\
  &\leq o\left(\frac{\delta^2}{T^{\frac{2+\nu}{2}}}\right)~,
\end{align}
where the last inequality is due to \eqref{si}.
\smallskip

\noindent Second, we have another estimation for the term $\sup_{ki}\left|\E \hat{\veps}_{ki}\right|$ as follows:
\begin{align}\label{bound2}
  &~\quad\sup_{ki}\left|\E \hat{\veps}_{ki}\right|=\sup_{ki}\Big|\E \big[\veps_{ki}\cdot \mathbf{1}_{\{|\veps_{ki}|\leq \delta T^{1/2}\}}\big]\Big|
  =\sup_{ki}\Big|\E \big[\veps_{ki}\cdot \mathbf{1}_{\{|\veps_{ki}|>\delta T^{1/2}\}}\big]\Big|\nonumber\\
  &\leq \frac{\sup_{ki}\E \big[|\veps_{ki}|^{4+\nu}\cdot \mathbf{1}_{\{|\veps_{ki}|>\delta T^{1/2}\}}\big]}{\left(\delta T^{1/2}\right)^{3+\nu}}=o\left(\frac{\delta}{T^{\frac{3+\nu}{2}}}\right)~.
\end{align}
Then similar to \eqref{j12}, we have
\begin{align*}
  &~\quad\left|\lambda_{\max}(\hat{B})-\lambda_{\max}(\tilde{B})\right|\nonumber\\
  &\leq\left\|\frac {1}{pT}\left(\hat{E}^T_1\hat{E}_1-\tilde{E}^T_1\tilde{E}_1\right)\hat{E}^T_2\hat{E}_2\right\|_{op}+\left\|\frac {1}{pT}\tilde{E}^T_1\tilde{E}_1\left(\hat{E}^T_2\hat{E}_2-\tilde{E}^T_2\tilde{E}_2\right)\right\|_{op}\nonumber\\
  &:=J_3+J_4~.
\end{align*}
Also, similar to \eqref{jjj} and \eqref{j11}, we have
\begin{align}
  &\quad\left\|\hat{E}^T_1\hat{E}_1-\tilde{E}^T_1\tilde{E}_1\right\|_{op}=\max_{\|x\|=1}x(\hat{E}^T_1\hat{E}_1-\tilde{E}^T_1\tilde{E}_1)x^T\nonumber\\
  &\leq\max_{\|x\|=1}x(\hat{E}^T_1\hat{E}_1-\tilde{E}^T_1\hat{E}_1)x^T+\max_{\|x\|=1}x(\tilde{E}^T_1\hat{E}_1-\tilde{E}^T_1\tilde{E}_1)x^T\nonumber\\
  &:=J_{31}+J_{32}~,
\end{align}
with
\begin{align}\label{wq}
  J_{31}&=\max_{\|x\|=1}x(\hat{E}^T_1\hat{E}_1-\tilde{E}^T_1\hat{E}_1)x^T
  =\max_{\|x\|=1}\sum_{i,j}x_ix_j\sum_{k=1}^p\left(\hat{\veps}_{ki}-\tilde{\veps}_{ki}\right)\hat{\veps}_{ki}\nonumber\\
  &\leq \bigg(\sum_{k=1}^p \sum_{i=1}^T\big(\hat{\veps}_{ki}-\tilde{\veps}_{ki}\big)^2\bigg)^{1/2}\cdot
  \bigg(\sum_{k=1}^p\sum_{j=1}^T\hat{\veps}^2_{kj}\bigg)^{1/2}\nonumber\\
  &=O\left(\left(pT\cdot \sum_{k=1}^p \sum_{i=1}^T\big(\hat{\veps}_{ki}-\tilde{\veps}_{ki}\big)^2\right)^{1/2}\right)
\end{align}
Since
\begin{align}
  &\quad~\sum_{k=1}^p \sum_{i=1}^T\big(\hat{\veps}_{ki}-\tilde{\veps}_{ki}\big)^2=\sum_{k=1}^p \sum_{i=1}^T\bigg(\hat{\veps}_{ki}-\frac{\hat{\veps}_{ki}-\E \hat{\veps}_{ki}}{\sigma_{ki}}\bigg)^2\nonumber\\
  &=\sum_{k=1}^p \sum_{i=1}^T\bigg(1-\frac{1}{\sigma_{ki}}\bigg)^2\hat{\veps}^2_{ki}+\sum_{k=1}^p \sum_{i=1}^T\frac{1}{\sigma^2_{ki}}\Big(\E\hat{\veps}_{ki}\Big)^2+\sum_{k=1}^p \sum_{i=1}^T\frac{2}{\sigma_{ki}}\Big(1-\frac{1}{\sigma_{ki}}\Big)\hat{\veps}_{ki}\E \hat{\veps}_{ki}\nonumber\\
  &\leq \max\left\{O\left(pT \cdot\left(\sup_{ki}\bigg|1-\frac{1}{\sigma_{ki}}\bigg|\right)^2\right), ~ O\left(pT\cdot\left(\sup_{ki}\Big|\E\hat{\veps}_{ki}\Big|\right)^2\right), \right.\nonumber\\
  &\quad \quad \quad \quad\left.O\left(pT\cdot\sup_{ki}\Big|1-\frac{1}{\sigma_{ki}}\Big|\cdot\sup_{ki}\big|\E \hat{\veps}_{ki}\big|\right)\right\}\nonumber\\
  &\leq\max \left\{ o\left(\frac{\delta^4 p}{T^{1+\nu}}\right),~o\left(\frac{\delta^2 p}{T^{2+\nu}}\right),~o\left(\frac{\delta^3 p}{T^{3/2+\nu}}\right)\right\}~,
\end{align}
where the last inequality is due to \eqref{bound1} and \eqref{bound2}.
Then according to \eqref{wq}, we have the bound for the term $J_{31}$:
\begin{align}\label{bound31}
  |J_{31}|\leq \max \left\{ o\left(\frac{\delta^2 p}{T^{\nu/2}}\right),~o\left(\frac{\delta p}{T^{\frac{1+\nu}{2}}}\right),~o\left(\frac{\delta\sqrt{\delta} p}{T^{1/4+\nu/2}}\right)\right\}~.
\end{align}
For the same reason, we have the term $|J_{32}|$ can be bounded by   \eqref{bound31} as well.

\noindent Therefore, we have
\begin{align*}
  |J_{3}|&=\left\|\frac {1}{pT}\left(\hat{E}^T_1\hat{E}_1-\tilde{E}^T_1\tilde{E}_1\right)\hat{E}^T_2\hat{E}_2\right\|_{op}\\
  &\leq \frac {1}{p}\left\|\hat{E}^T_1\hat{E}_1-\tilde{E}^T_1\tilde{E}_1\right\|_{op}\cdot \frac 1T\left\|\hat{E}^T_2\hat{E}_2\right\|_{op}\\
  &=O\left(\frac 1p(J_{31}+J_{32})\right)\\
  &\leq \max \left\{ o\left(\frac{\delta^2 }{T^{\nu/2}}\right),~o\left(\frac{\delta }{T^{\frac{1+\nu}{2}}}\right),~o\left(\frac{\delta\sqrt{\delta} }{T^{1/4+\nu/2}}\right)\right\}\rightarrow 0~.
\end{align*}

\noindent Similar, we also have $|J_4|\rightarrow 0$, which leads to the fact that
\begin{align}
  \left|\lambda_{\max}(\hat{B})-\lambda_{\max}(\tilde{B})\right|\rightarrow 0~.
\end{align}


\end{document}